\documentclass[12pt]{amsart}
\usepackage{amsmath, amssymb,graphicx,caption, subcaption, mathtools}

\DeclarePairedDelimiter\norm{\lVert}{\rVert}%

\theoremstyle{plain}

\newtheorem{theorem}{Theorem}[section]

\newtheorem{conjecture}{Conjecture}[section]
\theoremstyle{definition}

\theoremstyle{remark}
\newtheorem{remark}[theorem]{Remark}

\numberwithin{equation}{section}

\def\Dbar{\leavevmode\lower.6ex\hbox to 0pt{\hskip-.23ex \accent"16\hss}D}

\newcommand{\pp}{p}
\newcommand{\qq}{q}
\newcommand{\rr}{r}
\newcommand{\pe}{\tilde{p}}
\newcommand{\qe}{\tilde{q}}

\begin{document}

\title[Constructions of $SU(2)$ and Weyl equivariant maps]{Constructions of $SU(2)$ and Weyl equivariant maps for all classical groups}
%\titlerunning{On hyperbolic configurations of four points}
\author{Joseph Malkoun}
\address{Independent Researcher}

%\email{}


\maketitle

\begin{abstract} If $G$ is a compact Lie group, $T$ a maximal torus in $G$ (with Lie algebras $\mathfrak{g}$ and $\mathfrak{t}$ respectively) 
and $W$ the corresponding Weyl group, then the Berry-Robbins problem for $G$, as formulated by Sir Michael Atiyah and Roger Bielawski, asks whether there exists 
a continuous $SU(2) \times W$ equivariant map from the space of regular Cartan triples (an open subset of $\mathfrak{t} \otimes \mathbb{R}^3$) 
to $G/T$, where $SU(2)$ acts via a regular Lie group homomorphism $SU(2) \to G$. 
This was settled positively by Atiyah and Bielawski, and their maps are even smooth, 
but they are not explicit. For $G=U(n)$, there exists another construction due to Sir Michael Atiyah and developed further with Paul Sutcliffe, 
which is explicit, but relies on a linear independence conjecture. The author had previously found a similar type of construction 
for $G=Sp(m)$, also relying on a linear independence conjecture. In this paper, similar constructions are done for $SO(2m+1)$ and $SO(2m)$, 
thus exhausting the list of classical groups. \end{abstract}
\maketitle

\section{Introduction} \label{introduction}

The spin-statistics theorem says that for $n$ identical particles with spin $S$, where $S$ is a non-negative 
integer for bosons, and half a positive odd integer for fermions, then upon a complete interchange of the positions and spin states of any two of these $n$ particles, the 
wavefunction of the $n$ particles picks up a sign factor, $(-1)^{2S}$. In other words, in the case of $n$ identical bosons, interchanging 
the positions and spin states of two of the particles leaves the collective wavefunction invariant, while in the case of $n$ fermions, such an interchange produces 
a sign change in the collective wavefunction. In particular, this implies Fermi's exclusion principle: two fermions can not occupy simultaneously the same quantum state.

In \cite{BR1997}, M.V. Berry and J.M. Robbins were interested in a geometric proof of the spin-statistics theorem in quantum mechanics. Most of 
the ``standard'' proofs of that theorem rely on quantum field theory, while the theorem belongs to the realm of quantum mechanics. On the other hand, in their 
approach, Berry and Robbins explained the sign factor using a parallel transported spin basis, and the sign factor appears as a geometric phase, 
as one interchanges two particles via a smooth path in the configuration space. While generalizing from $2$ to $n$ particles, Berry and Robbins were led to ask the following question:

\emph{Does there exist for each $n \geq 2$, a continuous map 
\[ f_n: C_n(\mathbb{R}^3) \to U(n)/T^n, \] 
where $C_n(\mathbb{R}^3)$ is the configuration 
space of $n$ distinct points in $\mathbb{R}^3$, and $T^n$ is the $n$-torus of diagonal matrices in $U(n)$, which is equivariant for the 
action of the symmetric group $\Sigma_n$?}

A permutation $\sigma \in \Sigma_n$ acts on $\mathbf{x} = (\mathbf{x}_1, \, \ldots \, , \, \mathbf{x}_n) \in C_n(\mathbb{R}^3)$ as follows: 
\[ \sigma.(\mathbf{x}_1, \, \ldots \, , \, \mathbf{x}_n) := (\mathbf{x}_{\sigma^{-1}(1)}, \, \ldots \, , \, \mathbf{x}_{\sigma^{-1}(n)}) \]
and acts on $gT \in U(n)/T^n$, where $g\in U(n)$, by 
\[ \sigma.(gT) := gP_{\sigma}^{-1}T = g T P_{\sigma}^{-1}, \]
where $P_{\sigma}$ is the permutation matrix associated to $\sigma \in \Sigma_n$.

Recognizing $\Sigma_n$ as the Weyl group of $U(n)$, Atiyah and Bielawski asked the following generalization of the Berry-Robbins problem:

\emph{If $G$ is any compact Lie group, does there exist a continuous map $f_G: \mathfrak{t} \otimes \mathbb{R}^3 \setminus \Delta \to G/T$, 
where $\Delta$ is the union of the kernels of $\alpha \otimes 1_{\mathbb{R}^3}: \mathfrak{t} \otimes \mathbb{R}^3 \to \mathbb{R} \otimes \mathbb{R}^3 \simeq \mathbb{R}^3$, as $\alpha$ varies 
in the set of all roots of $\mathfrak{g}$, which is equivariant both for the action of the Weyl group $W$ and for that of $SU(2)$.}

Here, the Weyl group acts on $\mathfrak{t} \otimes \mathbb{R}^3 \setminus \Delta$ via its natural action on $\mathfrak{t}$, and its 
trivial action on $\mathbb{R}^3$, and acts on $G/T$ as follows: if $w \in N(T)$, then
\[ wT.(gT) = gw^{-1}T \]
which is well defined, as can be checked. On the other hand, $SU(2)$ acts on $\mathfrak{t} \otimes \mathbb{R}^3 \setminus \Delta$ via its trivial action on $\mathfrak{t}$ and its 
natural action on $\mathbb{R}^3$ (i.e. its adjoint action, which gives rise to the group homomorphism $SU(2) \to SO(3)$), and it acts on 
$G/T$ via a preferred regular homomorphism $\rho: SU(2) \to G$. A homomorphism from $SU(2)$ 
to $G$ is said to be \emph{regular} if its complexification takes a unipotent element in $SL(2,C)$ to a regular unipotent element in 
$G_{\mathbb{C}}$ (a unipotent element is said to be regular if it lies in a unique Borel subgroup). Such a homomorphism exists 
and is unique up to conjugation. Thus if $k \in SU(2)$, 
and $g \in G$, then
\[ k.gT := \rho(k)gT \]
Our task in this paper is to present smooth candidates of solutions of the generalized Berry-Robbins problem, in the sense of Atiyah 
and Bielawski, for the orthogonal groups $SO(2m+1)$ and $SO(2m)$. Such candidates are genuine solutions provided a linear independence conjecture holds, 
as in the original Atiyah-Sutcliffe construction (\cite{Atiyah2000}, \cite{Atiyah2001} and \cite{Atiyah-Sutcliffe2002}) 
which corresponds to $G = U(n)$. A similar construction was previously found by the author in \cite{Malkoun2014} for $G = Sp(m)$. Moreover, in \cite{Malkoun2019}, it was noticed 
that the linear independence conjectures for $U(2m)$ actually imply the corresponding linear independence conjectures for $Sp(m)$. 
We first review the known constructions for $G=U(n)$ and $G=Sp(m)$, before presenting the new constructions for 
$G=SO(2m+1)$ and $G=SO(2m)$.

\section{The unitary groups}

We review here the Atiyah-Sutcliffe construction (cf. \cite{Atiyah-Sutcliffe2002}), corresponding to $G=U(n)$. Let 
$\mathbf{x} = (\mathbf{x}_1, \, \ldots \, , \, \mathbf{x}_n) \in C_n(\mathbb{R}^3)$, where $C_n(\mathbb{R}^3)$ is the 
configuration space of $n$ distinct points in $\mathbb{R}^3$. In other words, $\mathbf{x}_1, \, \ldots \, , \, \mathbf{x}_n$ are 
$n$ distinct points in $\mathbb{R}^3$.

The Hopf map $h \colon S^3 \to S^2$, where $S^3$ is the unit sphere in $\mathbb{C}^2$, is defined by
\[ h(u, \, v) = ( 2 \, \Re(\bar{u} v), \, 2 \, \Im(\bar{u} v), \, |u|^2 - |v|^2 ). \]
Given $(x,\, y, \, z) \in S^2$, with $z \neq 1$, it can be checked that
\[ h^{-1}(x,\, y,\, z) = \frac{e^{i\theta}}{\sqrt{2(1-z)}} ( 1 - z, \, x+iy ), \text{ ($\theta \in \mathbb{R}$).} \]
The presence of the phase factor $e^{i \theta}$ is to be expected since $S^3$ is a (non-trivial) principal $U(1)$ bundle over $S^2$ with bundle 
map $h$. Moreover
\[ h^{-1}(0, \, 0, \, 1) = e^{i \theta} (0, \, 1) \text{ ($\theta \in \mathbb{R}$)}. \]
A choice of $(u,\, v) \in h^{-1}(x,\,y,\,z)$ is called a Hopf lift of $(x,\, y , \, z) \in S^2$.

For each $a$, $b$, $1 \leq a < b \leq n$, we form the vector $\omega_{ab} \in S^2$, obtained 
by ``looking'' from point $\mathbf{x}_a$ to point $\mathbf{x}_b$:
\[ \omega_{ab} = \frac{\mathbf{x}_b - \mathbf{x}_a}{\norm{\mathbf{x}_b-\mathbf{x}_a}} \]
Choose a Hopf lift $(u_{ab}, v_{ab}) \in S^3$ of $\omega_{ab}$ and then define the following complex polynomial depending on a complex variable $t$
\[ p_{ab}(t) = \left| \begin{array}{cc} u_{ab} & 1 \\ v_{ab} & t \end{array} \right| = u_{ab} \, t - v_{ab}. \]

We assume that once a choice of Hopf lift for $\omega_{ab}$ is made ($1 \leq a < b \leq n$), the following Hopf lift for $\omega_{ba} = - \omega_{ab}$ is selected
\[ (u_{ba}, \, v_{ba}) = (-\bar{v}_{ab}, \, \bar{u}_{ab}). \]
The reader may recognize that the right-hand side of the previous formula is obtained by applying the quaternionic structure $j$ to the chosen Hopf lift of $\omega_{ab}$. This step is important 
for defining the phase of the Atiyah-Sutcliffe determinant, as we will see shortly.

Form the following polynomials
\[ p_a(t) = \prod_{b \neq a} p_{ab}(t), \quad \text{($1 \leq a \leq n$)} \]
where the product is over all indices $b \in \{1, \, \ldots \, , \, n\}$ with $b \neq a$.

\begin{conjecture}[Atiyah-Sutcliffe conjecture $1$, cf. \cite{Atiyah2000},\cite{Atiyah2001} and \cite{Atiyah-Sutcliffe2002}] 
\label{conjectureA1} Given any 
configuration $\mathbf{x} \in C_n(\mathbb{R}^3)$, the $n$ polynomials $p_a(t)$, for $1 \leq a \leq n$, 
are $\mathbb{C}$-linearly independent.
\end{conjecture}
Provided this conjecture is true, the map $\mathbf{x} \to (p_1(t), \, \ldots \, , \, p_n(t))$ is a map from $C_n(\mathbb{R}^3)$ to 
$GL(n, \, \mathbb{C})/U(1)^n$, which is both $\Sigma_n$ and $SU(2)$ equivariant. Using polar decomposition, 
we can define a smooth map from $GL(n, \, \mathbb{C})/U(1)^n$ to $U(n)/T^n$, which is equivariant under the action of $\Sigma_n$ and 
$SU(2)$. More precisely, such a map is induced by the map $GL(n, \, \mathbb{C}) \to U(n)$ which maps
\[ g \mapsto (gg^*)^{-1/2} g \in U(n). \]

Thus, provided Conjecture \ref{conjectureA1} holds, one gets a smooth map 
\[ f_n: C_n(\mathbb{R}^3) \to U(n)/T^n, \]
which is both $\Sigma_n$ and $SU(2)$ equivariant.
Conjecture \ref{conjectureA1} was proved for $n=3$ by Atiyah (\cite{Atiyah2000} and \cite{Atiyah2001}) and for $n=4$ by 
Eastwood and Norbury in \cite{Eastwood-Norbury2001}, and by {\Dbar}okovi{\'c} for some special configurations 
(\cite{Dokovic2002a} and \cite{Dokovic2002b}). 

The Atiyah-Sutcliffe determinant (cf. \cite{Atiyah-Sutcliffe2002}) is defined as follows. Form the matrix
\[ M = (p_1, \, \ldots \, , \, p_n) \]
having as first column the coefficients of $p_1(t)$, ordered by increasing powers of $t$, and so on, and define
\[ D_{U(n)}(\mathbf{x}_1, \, \ldots \, , \, \mathbf{x}_n) := \det(M). \]

But we have made choices for half of the Hopf lifts (corresponding actually to a choice of positive roots of $U(n)$), so we need to check 
that $D_{U(n)}$ is independent of these choices. It can be checked, and we refer to \cite{Atiyah-Sutcliffe2002} for more details, that the phase of $D_{U(n)}$ is well defined 
precisely because of the fact that once Hopf lifts for the $\omega_{ab}$ are chosen, for $1 \leq a < b \leq n$, the Hopf lifts for $\omega_{ba}$ are then fixed to be 
$(-\bar{v}_{ab}, \, \bar{u}_{ab})$. So multiplying $(u_{ab}, v_{ab})$ by $e^{i\theta}$ results in multiplying $(u_{ba}, v_{ba})$ by $e^{-i \theta}$.

Moreover, $D_{U(n)}$ has the important properties of being invariant 
under the Weyl group $W$ of $G = U(n)$ ($W = \Sigma_n$ in this case), invariant 
under $SU(2)$, which acts on the configuration space via $SO(3)$, and invariant under the operation of scaling the original configuration $\mathbf{x}$ (cf. \cite{Atiyah-Sutcliffe2002}). It also 
gets conjugated by an odd orthogonal transformation in $\mathbb{R}^3$, such as a reflection with respect to a plane in $\mathbb{R}^3$, or the so called ``parity'' transformation $\mathbf{x} \mapsto -\mathbf{x}$.
 
The Atiyah-Sutcliffe conjecture $2$ can now be stated.
\begin{conjecture}[Atiyah-Sutcliffe conjecture $2$, \cite{Atiyah-Sutcliffe2002}] \label{conjectureA2} For any configuration 
$\mathbf{x} \in C_n(\mathbb{R}^3)$, we have $|D_{U(n)}(\mathbf{x})| \geq 1$. \end{conjecture}
This is a stronger conjecture than conjecture \ref{conjectureA1} (Atiyah-Sutcliffe conjecture $1$), and has been proved 
for $n = 3$ by Atiyah in \cite{Atiyah2000} and \cite{Atiyah2001}, and for $n=4$ by Bou Khuzam and Johnson 
in \cite{BKJ2014}, and by Svrtan independently and around the same time (cf. \cite{Svrtan2014}). At the time of writing this 
article, the Atiyah-Sutcliffe conjecture $2$ (Conjecture \ref{conjectureA2}) remains 
open for $n>4$, apart from some special configurations (cf. \cite{Dokovic2002a} and \cite{Dokovic2002b}).
\begin{remark} Though we shall not discuss it further, Atiyah and Sutcliffe made a third conjecture in 
\cite{Atiyah-Sutcliffe2002} (known as conjecture $3$), which, if true, would imply conjecture $2$. Conjecture $3$ was also proved to be 
true for $n=4$ in \cite{BKJ2014} and \cite{Svrtan2014}, independently and around the same time. \end{remark}

\section{The symplectic groups}

We review here the author's construction in \cite{Malkoun2014} for the Lie groups $G=Sp(m)$. If $\mathbf{x}_a \in \mathbb{R}^3$, 
for $1 \leq a \leq m$, then one can think of $\mathbf{x} = (\mathbf{x}_a)$ 
as an element of $\mathfrak{t} \otimes \mathbb{R}^3$, i.e. a triple of elements of $\mathfrak{t}$, where 
$\mathfrak{t} \simeq \mathbb{R}^m$ is the Lie algebra of a maximal torus $T^m$ of $G = Sp(m)$. The condition that $\mathbf{x}$ is a 
\emph{regular} Cartan triple amounts 
to requiring that $\mathbf{x}_a \neq 0$ for $1 \leq a \leq m$ and $\mathbf{x}_a \pm \mathbf{x}_b \neq \mathbf{0}$ 
for all $a$, $b$ with $1 \leq a < b \leq m$.

Choose for $1 \leq a < b \leq m$, the following Hopf lifts
\begin{align*} (u_{ab}^{-+}, \, v_{ab}^{-+}) & \in h^{-1}\left(\frac{-\mathbf{x}_a + \mathbf{x}_b}{\norm{-\mathbf{x}_a + \mathbf{x}_b}}\right) \\ 
(u_{ab}^{--}, \, v_{ab}^{--}) & \in h^{-1}\left(\frac{-\mathbf{x}_a - \mathbf{x}_b}{\norm{-\mathbf{x}_a - \mathbf{x}_b}}\right) \end{align*}

Once these choices are made, we then define
\begin{align*} (u_{ab}^{+-}, \, v_{ab}^{+-}) &= (-\bar{v}_{ab}^{-+}, \, \bar{u}_{ab}^{-+}) \\
(u_{ab}^{++}, \, v_{ab}^{++}) &= (-\bar{v}_{ab}^{--}, \, \bar{u}_{ab}^{--}) \end{align*}.

Also choose, for each $1 \leq a \leq n$,
\[ (u_a^{-}, \, v_a^{-}) \in h^{-1}\left(\frac{-\mathbf{x}_a}{\norm{-\mathbf{x}_a}}\right) \]
and then define
\[ (u_a^{+}, \, v_a^{+}) = (-\bar{v}_a^{-}, \, \bar{u}_a^{-}). \]

We then form, from each Hopf lift $(u, v)$, the corresponding linear polynomial $u\, t - v$, with the corresponding notation. 
We have thus defined $p_{ab}^{-+}(t)$, $p_{ab}^{--}(t)$, $p_{ab}^{+-}(t)$, $p_{ab}^{++}(t)$, for $1 \leq a < b \leq m$ and $p_a^{-}(t)$, $p_a^{+}(t)$ for $1 \leq a \leq m$.

We then extend the notation a little and also define, for $1 \leq a < b \leq m$, the following
\[ \begin{array}{cc} p_{ba}^{++}(t) = p_{ab}^{++}(t), & p_{ba}^{+-}(t) = p_{ab}^{-+}(t), \\
p_{ba}^{-+}(t) = p_{ab}^{+-}(t), & p_{ba}^{--}(t) = p_{ab}^{--}(t). \end{array} \]

We now define, for $a = 1, \, \ldots \, , \, m$, the polynomials $\pp_a(t)$ and $\qq_a(t)$ of degree at most $2m-1$.
\begin{align*} \pp_a(t) &= \prod_{b \neq a} p_{ab}^{-+}(t) \prod_{b \neq a} p_{ab}^{--}(t) \prod_{a = 1}^m p_a^{-}(t) \\
\qq_a(t) &= \prod_{b \neq a} p_{ab}^{+-}(t) \prod_{b \neq a} p_{ab}^{++}(t) \prod_{a = 1}^m p_a^{+}(t) \end{align*}

We then form the $2m \times 2m$ complex matrix $M$ as follows
\[ M = (\pp_1, \, \qq_1, \, \pp_2, \, \qq_2, \, \ldots, \, \pp_m, \,\qq_m) \]
which means for instance that the first column contains the polynomial $\pp_1(t)$, thought of as a $2m$-dimensional complex vector whose 
entries are the coefficients of $\pp_1(t)$, ordered by increasing powers of $t$, the second column 
contains the coefficients of $\qq_1(t)$, and so on.

\begin{conjecture}[Conjecture $1$ for $G = Sp(m)$, cf. \cite{Malkoun2014}] \label{conjectureB1} Given any configuration 
$\mathbf{x} \in (\mathfrak{t} \otimes \mathbb{R}^3) \setminus \Delta$, corresponding to $G = Sp(m)$ (in other 
words, $\mathbf{x}_a \neq 0$ for all $1 \leq a \leq m$, and $\mathbf{x}_a \pm \mathbf{x}_b \neq 0$, for all 
$1 \leq a < b \leq m$), the corresponding polynomials $\pp_1, \, \qq_1, \, \ldots \, , \, \pp_m, \, \qq_m$ are linearly independent 
over $\mathbb{C}$.
\end{conjecture}

If Conjecture $1$ for $G = Sp(m)$ is true (in other words, Conjecture \ref{conjectureB1}), then one can show that we have constructed a smooth map with domain 
the configuration space $(\mathfrak{t} \otimes \mathbb{R}^3) \setminus \Delta$ of regular Cartan triples and target 
$GL(m, \, \mathbb{H})/U(1)^m$, which is equivariant under the Weyl group $W$ of $G = Sp(m)$. This makes 
use of the observation that a pair $(\pp_a, \qq_a)$ defines a quaternionic 
vector $v_a \in (\mathbb{H}^m \setminus \{ \mathbf{0} \})/U(1)$, since the roots of $q_a$ are the antipodals 
of those of $p_a$, and the antipodal map, in the case of an odd number of roots, induces a quaternionic map on the 
corresponding polynomial space. For more details, the reader may refer to \cite{Malkoun2014}. 

We now construct a smooth $N(T^m)$-equivariant map $GL(m, \, \mathbb{H}) \to Sp(m)$, where $N(T^m)$ is the normalizer of the maximal torus $T^m$ in $Sp(m)$, given by
\[ g \mapsto (gg^*)^{-1/2} \, g \]
where $N(T^m)$ acts on both spaces by multiplication from the right and $g^*$ denotes now the \emph{quaternionic} conjugate transpose of $g \in GL(m, \, \mathbb{H})$. This map descends to a $W$-equivariant map 
\[ GL(m, \, \mathbb{H})/U(1)^m \to Sp(m)/T^m. \]

Thus, provided Conjecture \ref{conjectureB1} is true, the composition of the previous two maps
\[ (\mathfrak{t} \otimes \mathbb{R}^3) \setminus \Delta \to GL(m, \, \mathbb{H})/U(1)^m \to Sp(m)/T^m \]
gives the required $SU(2) \times W$ equivariant map, where $k \in SU(2)$ acts $Sp(m)/T^m$ by left multiplication by $\rho(k)$, where $\rho: SU(2) \to Sp(m)$ is a regular homomorphism. 

We shall also define a complex-valued normalized determinant function $D_{Sp(m)}$ on $(\mathfrak{t} \otimes \mathbb{R}^3) \setminus \Delta $ by
\[ D_{Sp(m)}(\mathbf{x}_1, \, \ldots \, , \, \mathbf{x}_m) = \det(M) = \det (\pp_1, \, \qq_1, \, \ldots \, , \, \pp_m, \,\qq_m).\]
The normalized determinant $D_{Sp(m)}$ is invariant under the Weyl group $W$ of $G = Sp(m)$, as well as the action of 
$SO(3)$ and the operation of scaling the configuration $\mathbf{x}$ (cf \cite{Malkoun2014} for proofs of these claims).

\begin{conjecture}[Conjecture $2$ for $G = Sp(m)$, \cite{Malkoun2014}] \label{conjectureB2} For any configuration 
$\mathbf{x} \in (\mathfrak{t} \otimes \mathbb{R}^3)\setminus \Delta$, for $G = Sp(m)$ (which means that 
$\mathbf{x}_a \neq \mathbf{0}$ for all $1 \leq a \leq m$, and $\mathbf{x}_a \pm \mathbf{x}_b \neq \mathbf{0}$, for 
all $1 \leq a < b \leq m$), we have $|D_{Sp(m)}(\mathbf{x})| \geq 1$. \end{conjecture}
This was actually proved for $m=2$ in \cite{Malkoun2014}, and is a conjecture for $m \geq 3$ (at least at the time 
of writing).

Moreover, in \cite{Malkoun2019}, it was shown that if $D_{U(2m)}$ is non-vanishing on $C_{2m}(\mathbb{R}^3)$, then $D_{Sp(m)}$ is non-vanishing on its domain. It was also shown that if $|D_{U(2m)}| \geq 1$ for any 
$\mathbf{x} \in C_{2m}(\mathbb{R}^3)$, then $|D_{Sp(m)}| \geq 1$ on its domain. So the linear independence conjectures for $U(2m)$ actually imply the corresponding linear independence conjectures for $Sp(m)$. This follows 
from an interesting and elementary connection between the roots of $U(2m)$ and those of $Sp(m)$ (cf. \cite{Malkoun2019}).        

\section{The orthogonal groups in odd dimensions}

We present here a new construction (though known to the author for several years actually) for $G = SO(2m+1)$. If $\mathbf{x}_a \in \mathbb{R}^3$, 
for $1 \leq a \leq m$, then one can think of $\mathbf{x} = (\mathbf{x}_a)$ 
as an element of $\mathfrak{t} \otimes \mathbb{R}^3$, i.e. a triple of elements of $\mathfrak{t}$, where 
$\mathfrak{t} \simeq \mathbb{R}^m$ is the Lie algebra of a maximal torus $T^m$ of $G = SO(2m+1)$. The condition that $\mathbf{x}$ is a 
\emph{regular} Cartan triple amounts 
to requiring that $\mathbf{x}_a \neq 0$ for $1 \leq a \leq m$ and $\mathbf{x}_a \pm \mathbf{x}_b \neq \mathbf{0}$ 
for all $a$, $b$ with $1 \leq a < b \leq m$.

As in the previous section, choose for $1 \leq a < b \leq m$, the following Hopf lifts
\begin{align*} (u_{ab}^{-+}, \, v_{ab}^{-+}) & \in h^{-1}\left(\frac{-\mathbf{x}_a + \mathbf{x}_b}{\norm{-\mathbf{x}_a + \mathbf{x}_b}}\right) \\ 
(u_{ab}^{--}, \, v_{ab}^{--}) & \in h^{-1}\left(\frac{-\mathbf{x}_a - \mathbf{x}_b}{\norm{-\mathbf{x}_a - \mathbf{x}_b}}\right) \end{align*}

Once these choices are made, we then define
\begin{align*} (u_{ab}^{+-}, \, v_{ab}^{+-}) &= (-\bar{v}_{ab}^{-+}, \, \bar{u}_{ab}^{-+}) \\
(u_{ab}^{++}, \, v_{ab}^{++}) &= (-\bar{v}_{ab}^{--}, \, \bar{u}_{ab}^{--}) \end{align*}.

Also choose, for each $1 \leq a \leq n$,
\[ (u_a^{-}, \, v_a^{-}) \in h^{-1}\left(\frac{-\mathbf{x}_a}{\norm{-\mathbf{x}_a}}\right) \]
and then define
\[ (u_a^{+}, \, v_a^{+}) = (-\bar{v}_a^{-}, \, \bar{u}_a^{-}). \]

We then form, from each Hopf lift $(u, v)$, the corresponding linear polynomial $u\, t - v$, with the corresponding notation. 
We have thus defined $p_{ab}^{-+}(t)$, $p_{ab}^{--}(t)$, $p_{ab}^{+-}(t)$, $p_{ab}^{++}(t)$, for $1 \leq a < b \leq m$ and $p_a^{-}(t)$, $p_a^{+}(t)$ for $1 \leq a \leq m$.

We then extend the notation a little and also define, for $1 \leq a < b \leq m$, the following
\[ \begin{array}{cc} p_{ba}^{++}(t) = p_{ab}^{++}(t), & p_{ba}^{+-}(t) = p_{ab}^{-+}(t), \\
p_{ba}^{-+}(t) = p_{ab}^{+-}(t), & p_{ba}^{--}(t) = p_{ab}^{--}(t). \end{array} \]

We now define, for $a = 1, \, \ldots \, , \, m$, the polynomials $\pp_a(t)$ and $\qq_a(t)$ of degree at most $2m$.
\begin{align*} \pp_a(t) &= \prod_{b \neq a} p_{ab}^{-+}(t) \prod_{b \neq a} p_{ab}^{--}(t) \prod_{a = 1}^m p_a^{-}(t)^2 \\
\qq_a(t) &= \prod_{b \neq a} p_{ab}^{+-}(t) \prod_{b \neq a} p_{ab}^{++}(t) \prod_{a = 1}^m p_a^{+}(t)^2 .\end{align*}

We need one more polynomial, also of degree at most $2m$, which we now define.
\[ \rr(t) = \prod_{a=1}^m p_a^{-}(t) \prod_{a=1}^m p_a^{+}(t). \]

We then form the $(2m+1) \times (2m+1)$ complex matrix $M$ as follows
\[ M = (\pp_1, \, \qq_1, \, \pp_2, \, \qq_2, \, \ldots, \, \pp_m, \,\qq_m, \, \rr) \]
which means for instance that the first column contains the polynomial $\pp_1(t)$, thought of as a $2m+1$-dimensional complex vector whose 
entries are the coefficients of $\pp_1(t)$, ordered by increasing powers of $t$, the second column 
contains the coefficients of $\qq_1(t)$, and so on.

\begin{conjecture}[Conjecture $1$ for $G = SO(2m+1)$] \label{conjectureC1} Given any configuration 
$\mathbf{x} \in (\mathfrak{t} \otimes \mathbb{R}^3) \setminus \Delta$, corresponding to $G = SO(2m+1)$ (in other 
words, $\mathbf{x}_a \neq 0$ for all $1 \leq a \leq m$, and $\mathbf{x}_a \pm \mathbf{x}_b \neq 0$, for all 
$1 \leq a < b \leq m$), the corresponding polynomials $\pp_1, \, \qq_1, \, \ldots \, , \, \pp_m, \, \qq_m$ and $\rr$ are linearly independent 
over $\mathbb{C}$.
\end{conjecture}

If Conjecture $1$ for $G = SO(2m+1)$ is true (in other words, Conjecture \ref{conjectureC1}), then one can show that we have constructed a smooth map with domain 
the configuration space $(\mathfrak{t} \otimes \mathbb{R}^3) \setminus \Delta$ of regular Cartan triples and target 
$GL(2m+1, \, \mathbb{R})/SO(2)^m$, which is equivariant under the Weyl group $W$ of $G = SO(2m+1)$. This makes 
use of the observation that a linearly independent (over $\mathbb{C}$) pair $(\pp_a, \qq_a)$ defines a complex $2$-dimensional subspace $V_a$ of $\mathbb{C}^{2m+1}$, which is 
preserved by the real structure $\tau$ on $\mathbb{C}^{2m+1}$, itself induced by the antipodal map on the set of roots (taking into account multiplicity). Indeed, the number of roots of $\pp_a$ and 
$\qq_a$ is $2m$ which is even, so the antipodal map induces a real structure on the corresponding polynomial space (if one has a quaternionic structure, say $j$, on $\mathbb{C}^2$, then the 
symmetric tensor of $2m$ copies of $j$ is then a real structure on the $2m$-th symmetric tensor power $S^{2m}(\mathbb{C}^2)$ of $\mathbb{C}^2$). Thus the intersection of $V_a$ with the real slice of 
$\mathbb{C}^{2m+1}$ (with respect to the real structure $\tau$) gives a real $2$-dimensional subspace of $\mathbb{R}^{2m+1}$. Note also that $\rr$ is preserved by $\tau$, so it is a \emph{real} 
vector with respect to $\tau$.

We now construct a smooth $N(T^m)$-equivariant map 
\[ GL^+(2m+1, \, \mathbb{R}) \to SO(2m+1), \] 
where $N(T^m)$ is the normalizer of the maximal torus $T^m$ in $SO(2m+1)$, given by
\[ g \mapsto (gg^T)^{-1/2} \, g \]
where $N(T^m)$ acts on both spaces by multiplication from the right and $g^T$ denotes the transpose of $g \in GL^+(2m+1, \, \mathbb{R})$. This map descends to a $W$-equivariant map 
\[ GL^+(2m+1, \, \mathbb{R})/SO(2)^m \to SO(2m+1)/T^m. \]

Conjecture \ref{conjectureB1} only gives us a map into $GL(2m+1, \mathbb{R})/SO(2)^m$, which has $2$ connected components, corresponding to positive/negative determinant 
matrices. However, the domain of that map (the space of regular Cartan triples) is itself connected, so the image of our map would then be contained in exactly one of these 
two connected components. In the case where say (hypothetically), the images under our maps 
all have negative determinants, then multiplying the images by a sign factor would then make the new images have positive determinants, while preserving $SU(2) \times W$ equivariance. 

Thus, provided Conjecture \ref{conjectureB1} is true, the composition of the previous two maps
\[ (\mathfrak{t} \otimes \mathbb{R}^3) \setminus \Delta \to GL^+(2m+1, \, \mathbb{R})/U(1)^m \to SO(2m+1)/T^m \]
gives the required $SU(2) \times W$ equivariant map, where $k \in SU(2)$ acts $SO(2m+1)/T^m$ by left multiplication by $\rho(k)$, where $\rho: SU(2) \to SO(2m+1)$ is a regular homomorphism. Essentially, 
the $SU(2)$ equivariance property of our map stems from the fact that the Hopf map is itself $SU(2)$ equivariant, for the natural action of $SU(2)$ on $\mathbb{C}^2$ and its adjoint representation on $\mathbb{R}^3$.

We shall also define a normalized determinant function 
\[ D_{SO(2m+1)} \colon (\mathfrak{t} \otimes \mathbb{R}^3) \setminus \Delta \to \mathbb{C} \] 
by
\[ D_{SO(2m+1)}(\mathbf{x}_1, \, \ldots \, , \, \mathbf{x}_m) = \det(M) = \det (\pp_1, \, \qq_1, \, \ldots \, , \, \pp_m, \, \qq_m, \, \rr).\]
We claim that the normalized determinant $D_{SO(2m+1)}$ is invariant under the Weyl group $W$ of $G = SO(2m+1)$. Indeed, if we permute two of the points of a regular Cartan triple $\mathbf{x}$ , say the first two points, 
$\mathbf{x}_1$ and $\mathbf{x}_2$, then the determinant of $M$ picks up $2$ minus signs (due to interchanging $\pp_1$, $\qq_1$ with $\pp_2$, $\qq_2$ respectively), but there are 
other minus signs, which are more subtle, coming from our convention that 
\[ (u^{+-}_{12}, \, v^{+-}_{12}) = (- \bar{v}^{-+}_{12}, \, \bar{u}^{-+}_{12}). \]
The latter minus signs are due to the fact that the quaternionic structure $j$ on $\mathbb{C}^2$ (corresponding to the antipodal map on $S^2$) satisfies $j^2 = -1_{\mathbb{C}^2}$. There are actually $2$ such extra 
minus signs, since $p^{+-}_{12}(t)$ occurs once in $\pp_1$ and once in $\qq_2$. Hence $D_{SO(2m+1)}(\mathbf{x})$ is invariant under 
the operation of interchanging two of the points in $\mathbf{x} \in (\mathfrak{t} \otimes \mathbb{R}^3) \setminus \Delta$, and so it is invariant under the action of the permutation subgroup $\Sigma_m \subset W$.

Similarly, if we replace $\mathbf{x}_1$ with $-\mathbf{x}_1$, we then get a minus sign in the determinant of $M$ from interchanging $\pp_1$ and $\qq_1$. But there are an odd number of extra minus signs, 
essentially due to the fact that $p^{+}_1(t)$ appears twice as a factor of $\pp_1(t)$ and once as a factor of $\rr(t)$, so that it appears an odd number of times.

So we have thus showed that $D_{SO(2m+1)}$ is invariant under the action of the Weyl group $W$. It can also be seen that $D_{SO(2m+1)}$ is invariant under the action of any $k \in SO(3)$, as $M$ would then become 
$\rho(\tilde{k}) M$, where $\tilde{k} \in SU(2)$ is a``lift'' of $k$ (since $SU(2)$ is a double cover of $SO(3)$) and $\rho: SU(2) \to SO(2m+1)$ is a regular homomorphism, which actually factors through $SO(3)$ in this case. Thus 
$D_{SO(2m+1)}$ is invariant under the action of $SO(3)$, since $\operatorname{det}(\rho(\tilde{k})) = 1$.

Similarly to the other cases, we make the following conjecture.

\begin{conjecture}[Conjecture $2$ for $G = SO(2m+1)$] \label{conjectureC2} For any configuration 
$\mathbf{x} \in (\mathfrak{t} \otimes \mathbb{R}^3)\setminus \Delta$, for $G = SO(2m+1)$ (which means that 
$\mathbf{x}_a \neq \mathbf{0}$ for all $1 \leq a \leq m$, and $\mathbf{x}_a \pm \mathbf{x}_b \neq \mathbf{0}$, for 
all $1 \leq a < b \leq m$), we have $|D_{SO(2m+1)}(\mathbf{x})| \geq 1$. \end{conjecture}

We have verified this conjecture numerically for small values of $m$, by generating each time a pseudo-random sample of $1000$ points in $(\mathfrak{t} \otimes \mathbb{R}^3) \setminus \Delta$ and checking that 
the inequality $|D_{SO(2m+1)}(\mathbf{x})| \geq 1$ indeed holds for that sample.

\section{The orthogonal groups in even dimensions}

We present here a new construction for $G = SO(2m)$. It is in the same spirit as the other constructions, but a little more involved. In fact, the author discovered it only recently (at the time of 
writing) and several years after discovering the construction in the previous section for $SO(2m+1)$ (though he had a preliminary unsatisfactory construction for $SO(2m)$ before). If $\mathbf{x}_a \in \mathbb{R}^3$, 
for $1 \leq a \leq m$, then one can think of $\mathbf{x} = (\mathbf{x}_a)$ 
as an element of $\mathfrak{t} \otimes \mathbb{R}^3$, i.e. a triple of elements of $\mathfrak{t}$, where 
$\mathfrak{t} \simeq \mathbb{R}^m$ is the Lie algebra of a maximal torus $T^m$ of $G = SO(2m)$. The condition that $\mathbf{x}$ is a 
\emph{regular} Cartan triple amounts to requiring that $\mathbf{x}_a \pm \mathbf{x}_b \neq \mathbf{0}$ for all $a$, $b$ with $1 \leq a < b \leq m$.

As in the previous section, choose for $1 \leq a < b \leq m$, the following Hopf lifts
\begin{align*} (u_{ab}^{-+}, \, v_{ab}^{-+}) & \in h^{-1}\left(\frac{-\mathbf{x}_a + \mathbf{x}_b}{\norm{-\mathbf{x}_a + \mathbf{x}_b}}\right) \\ 
(u_{ab}^{--}, \, v_{ab}^{--}) & \in h^{-1}\left(\frac{-\mathbf{x}_a - \mathbf{x}_b}{\norm{-\mathbf{x}_a - \mathbf{x}_b}}\right) \end{align*}

Once these choices are made, we then define
\begin{align*} (u_{ab}^{+-}, \, v_{ab}^{+-}) &= (-\bar{v}_{ab}^{-+}, \, \bar{u}_{ab}^{-+}) \\
(u_{ab}^{++}, \, v_{ab}^{++}) &= (-\bar{v}_{ab}^{--}, \, \bar{u}_{ab}^{--}) \end{align*}.

We then form, from each Hopf lift $(u, v)$, the corresponding linear polynomial $u\, t - v$, with the corresponding notation. 
We have thus defined $p_{ab}^{-+}(t)$, $p_{ab}^{--}(t)$, $p_{ab}^{+-}(t)$ and $p_{ab}^{++}(t)$, for $1 \leq a < b \leq m$.

We then extend the notation a little and also define, for $1 \leq a < b \leq m$, the following
\[ \begin{array}{cc} p_{ba}^{++}(t) = p_{ab}^{++}(t), & p_{ba}^{+-}(t) = p_{ab}^{-+}(t), \\
p_{ba}^{-+}(t) = p_{ab}^{+-}(t), & p_{ba}^{--}(t) = p_{ab}^{--}(t). \end{array} \]

We now define, for $a = 1, \, \ldots \, , \, m$, the polynomials $\pp_a(t)$ and $\qq_a(t)$ of degree at most $2m-2$.
\begin{align*} \pp_a(t) &= \prod_{b \neq a} p_{ab}^{-+}(t) \prod_{b \neq a} p_{ab}^{--}(t) \\
\qq_a(t) &= \prod_{b \neq a} p_{ab}^{+-}(t) \prod_{b \neq a} p_{ab}^{++}(t) .\end{align*}

We then form the $(2m-1) \times 2m$ complex matrix $M'$ as follows
\[ M' = (\pp_1, \, \qq_1, \, \pp_2, \, \qq_2, \, \ldots, \, \pp_m, \,\qq_m) \]
which means for instance that the first column contains the polynomial $\pp_1(t)$, thought of as a $2m-1$-dimensional complex vector whose 
entries are the coefficients of $\pp_1(t)$, ordered by increasing powers of $t$, the second column contains the coefficients of $\qq_1(t)$, and so on.

Note that $M'$ is not a square matrix, but we will construct a new $2m \times 2m$ matrix $M$ containing $M'$ as its upper $(2m-1) \times 2m$ block. It remains to define the last row 
of $M$.

Equivalently, for $1 \leq a \leq m$, we could extend each $\pp_a$, thought of as a complex vector with $2m-1$ coordinates (the coefficients of the polynomial $\pp_a(t)$), to a complex vector denoted by 
$\pe_a$, whose first $2m-1$ coordinates are the coefficients of $\pp_a$. Let us denote the last coordinate of $\pe_a$ by $c^{(a)}$. Similarly, we extend each $\qq_a$ to a $2m$-dimensional complex 
vector, which we denote by $\qe_a$, and we denote its last coordinate by $d^{(a)}$. We can now define, for $a = 1, \, \ldots \, , \, m$,
\begin{align*} c^{(a)} & = -\frac{1}{2} \, \prod_{b \neq a} \operatorname{det}(p^{-+}_{ab}, \, p^{--}_{ab}) \\
d^{(a)} & = \frac{1}{2} \, \prod_{b \neq a} \operatorname{det}(p^{++}_{ab}, \, p^{+-}_{ab}) \end{align*}

\begin{conjecture}[Conjecture $1$ for $G = SO(2m)$] \label{conjectureD1} Given any configuration 
$\mathbf{x} \in (\mathfrak{t} \otimes \mathbb{R}^3) \setminus \Delta$, corresponding to $G = SO(2m)$ (in other 
words, $\mathbf{x}_a \pm \mathbf{x}_b \neq 0$, for all $1 \leq a < b \leq m$), the corresponding vectors $\pe_1, \, \qe_1, \, \ldots \, , \, \pe_m, \, \qe_m$ (each of these vectors only defined 
up to a phase factor) are linearly independent over $\mathbb{C}$.
\end{conjecture}

If Conjecture $1$ for $G = SO(2m)$ is true (in other words, Conjecture \ref{conjectureD1}), then one can show that we have constructed a smooth map with domain 
the configuration space $(\mathfrak{t} \otimes \mathbb{R}^3) \setminus \Delta$ of regular Cartan triples and target $GL(2m, \, \mathbb{R})/SO(2)^m$, which is equivariant under the 
Weyl group $W$ of $G = SO(2m)$. This makes use of the observation that a linearly independent (over $\mathbb{C}$) pair $(\pe_a, \qe_a)$ defines a complex $2$-dimensional subspace $V_a$ of $\mathbb{C}^{2m}$, which is 
preserved by the real structure $\tau$ on $\mathbb{C}^{2m}$ defined as follows. The real structure $\tau$ on $\mathbb{C}^{2m}$ is induced by the antipodal map on the set of roots (taking into account multiplicity) on the polynomial space (of polynomials 
of degree at most $2m-2$) corresponding to the first $2m-1$ coordinates, and which acts by mapping the last coordinate (i.e. the $2m$-th coordinate), say $z_{2m}$, to $(-1)^m \bar{z}_{2m}$. 
Thus the intersection of $V_a$ with the real slice of $\mathbb{C}^{2m}$ (with respect to the real structure $\tau$) gives a real $2$-dimensional subspace of $\mathbb{R}^{2m}$.

We now construct a smooth $N(T^m)$-equivariant map 
\[ GL^+(2m, \, \mathbb{R}) \to SO(2m), \]
where $N(T^m)$ is the normalizer of the maximal torus $T^m$ in $SO(2m)$, given by
\[ g \mapsto (gg^T)^{-1/2} \, g \]
where $N(T^m)$ acts on both spaces by multiplication from the right and $g^T$ denotes the transpose of $g \in GL^+(2m, \, \mathbb{R})$. This map descends to a $W$-equivariant map 
\[ GL^+(2m, \, \mathbb{R})/SO(2)^m \to SO(2m)/T^m. \]

Conjecture \ref{conjectureD1} only gives us a map into $GL(2m, \mathbb{R})/SO(2)^m$, which has $2$ connected components, corresponding to posivite/negative determinant 
matrices. However, just as in the previous section, the domain of that map (the space of regular Cartan triples) is itself connected, so the image of our map would then be contained in exactly one of these 
two connected components. In the case where say (hypothetically), the images under our maps 
all have negative determinants, then multiplying the last row of $M$ by a sign factor would then make the new images have positive determinants, while preserving $SU(2) \times W$ equivariance. 

Thus, provided Conjecture \ref{conjectureD1} is true, the composition of the previous two maps
\[ (\mathfrak{t} \otimes \mathbb{R}^3) \setminus \Delta \to GL^+(2m, \, \mathbb{R})/U(1)^m \to SO(2m)/T^m \]
gives the required $SU(2) \times W$ equivariant map, where $k \in SU(2)$ acts $SO(2m)/T^m$ by left multiplication by $\rho(k)$, where $\rho: SU(2) \to SO(2m)$ is a regular homomorphism. 

We shall also define a normalized determinant function 
\[ D_{SO(2m)} \colon (\mathfrak{t} \otimes \mathbb{R}^3) \setminus \Delta \to \mathbb{C} \] 
by
\[ D_{SO(2m)}(\mathbf{x}_1, \, \ldots \, , \, \mathbf{x}_m) = \det(M) = \det (\pe_1, \, \qe_1,  \, \ldots, \, \pe_m, \, \qe_m).\]
We claim that the normalized determinant $D_{SO(2m)}$ is invariant under the Weyl group $W$ of $G = SO(2m)$. Indeed, if we permute two of the points of a regular Cartan triple $\mathbf{x}$ , say the first two points, 
$\mathbf{x}_1$ and $\mathbf{x}_2$, then the determinant of $M$ picks up $2$ minus signs (due to interchanging $\pe_1$, $\qe_1$ with $\pe_2$, $\qe_2$ respectively), but there are 
other minus signs, which are more subtle, coming from our convention that 
\[ (u^{+-}_{12}, \, v^{+-}_{12}) = (- \bar{v}^{-+}_{12}, \, \bar{u}^{-+}_{12}). \]
The latter minus signs are due to the fact that the quaternionic structure $j$ on $\mathbb{C}^2$ (corresponding to the antipodal map on $S^2$) satisfies $j^2 = -1_{\mathbb{C}^2}$. There are actually $2$ such extra 
minus signs, coming from the fact that $p^{+-}_{12}(t)$ occurs once in $\pp_1$ and once in $\qq_2$. Hence $D_{SO(2m)}(\mathbf{x})$ is invariant under 
transpositions of two of the points in $\mathbf{x} \in (\mathfrak{t} \otimes \mathbb{R}^3) \setminus \Delta$, and so it is invariant under the action of the permutation subgroup $\Sigma_m \subset W$.

Similarly, if we simultaneously replace $\mathbf{x}_1$, $\mathbf{x}_2$ with $-\mathbf{x}_1$ and $-\mathbf{x}_2$ respectively, we then get $2$ minus signs in the determinant of $M$ 
from interchanging the columns $\pe_1$, $\pe_2$ with $\qe_1$ and $\qe_2$ respectively. And there is an even number of extra minus signs, since for instance $p^{-+}_{12}(t)$ occurs once in $\pe_1$ and once in 
$\qe_2$ (so an even number of times), and so on.

So we have thus showed that $D_{SO(2m)}$ is invariant under the action of the Weyl group $W$. It can also be seen that $D_{SO(2m)}$ is invariant under the action of any $k \in SO(3)$, as $M$ would then become 
$\rho(\tilde{k}) M$, where $\tilde{k} \in SU(2)$ is a ``lift'' of $k$ (since $SU(2)$ is a double cover of $SO(3)$) and $\rho: SU(2) \to SO(2m)$ is a regular homomorphism, which actually factors through $SO(3)$ in this case. 
We are also implicitly making use of the fact that the $c^{(a)}$ and $d^{(a)}$ are invariant under the action of $SO(3)$, since they are products of $2$ by $2$ determinants built from the linear factors $p^{-+}_{ab}$, 
$p^{--}_{ab}$ and so on. Thus $D_{SO(2m)}$ is invariant under the action of $SO(3)$, since $\operatorname{det}(\rho(\tilde{k})) = 1$.

Similarly to the other cases, we make the following conjecture.

\begin{conjecture}[Conjecture $2$ for $G = SO(2m)$] \label{conjectureD2} For any configuration 
$\mathbf{x} \in (\mathfrak{t} \otimes \mathbb{R}^3)\setminus \Delta$, for $G = SO(2m)$ (which means that 
$\mathbf{x}_a \pm \mathbf{x}_b \neq \mathbf{0}$, for all $1 \leq a < b \leq m$), we have $|D_{SO(2m)}(\mathbf{x})| \geq 1$. \end{conjecture}

We have also verified this conjecture numerically for small values of $m$, by generating each time a pseudo-random sample of $1000$ points in $(\mathfrak{t} \otimes \mathbb{R}^3) \setminus \Delta$ and checking that 
the inequality $|D_{SO(2m)}(\mathbf{x})| \geq 1$ indeed holds for that sample.

\section{Concluding Remarks}

In order to define our maps, we actually made use of ideas contained in \cite{Malkoun2020}, though the matrices one obtains there are usually not square, so the author had to work harder in order to define the maps here, 
particularly so for the maps corresponding to $SO(2m)$ (which are a little more involved than for the other cases). Indeed, one had to add an extra coordinate to each column, and it was not clear 
at all how to do so. The author made use of Python $3$ programs in order to check the conjectures numerically. The maps constructed here are, in the author's opinion, natural, partly since they satisfy $SU(2) \times W$ equivariance. They are also explicit.

In the wonderful article \cite{Atiyah-Bielawski2002}, the authors identify the space of regular Cartan triples essentially as a space of solutions of Nahm's equations on the half-line $(0, \, \infty)$ satisfying some specific boundary conditions. They 
then obtain a solution of the Berry-Robbins problem by mapping a solution of Nahm's equations to its asymptotics at infinity. It would be interesting to investigate whether or not their solutions and our constructions are related in some way. This question 
was actually asked by Sir Michael Atiyah for $G = U(n)$ in his Edinburgh Lectures on Geometry, Analysis and Physics (Problem $1$ towards the end of section $1.8$ in \cite{Atiyah08-09}).

\providecommand{\bysame}{\leavevmode\hbox to3em{\hrulefill}\thinspace}
\providecommand{\MR}{\relax\ifhmode\unskip\space\fi MR }
\providecommand{\MRhref}[2]{%
  \href{http://www.ams.org/mathscinet-getitem?mr=#1}{#2}
}
\providecommand{\href}[2]{#2}


\begin{thebibliography}{1}

\bibitem{Atiyah2000}
Atiyah, M.F. \emph{The geometry of classical particles}, Surveys in
  differential geometry, Surv. Differ. Geom., VII, Int. Press, Somerville, MA,
  2000, pp.~1--15.

\bibitem{Atiyah2001}
Atiyah, M.F. \emph{Configurations of points}, R. Soc. Lond. Philos. Trans. Ser. A
  Math. Phys. Eng. Sci. \textbf{359} (2001), no.~1784, 1375--1387, Topological
  methods in the physical sciences (London, 2000). %\MR{1853626 (2003b:55016)}
  
\bibitem{Atiyah08-09} \emph{Edinburgh Lectures on Geometry, Analysis and Physics}, (math) arXiv: 1009.4827.

\bibitem{Atiyah-Bielawski2002}
Atiyah, M.F. and Bielawski, R. \emph{Nahm's equations, configuration spaces and flag manifolds}, 
Bull. Braz. Math. Soc. (N.S.) \textbf{33} (2002), no. 2, 157--176.

\bibitem{Atiyah-Sutcliffe2002}
Atiyah, M.F. and Sutcliffe, P.M. \emph{The geometry of point particles}, R.
  Soc. Lond. Proc. Ser. A Math. Phys. Eng. Sci. \textbf{458} (2002), no.~2021,
  1089--1115. %\MR{1902577 (2003c:55019)}

\bibitem{BR1997}
Berry, M.V. and Robbins, J.M. \emph{Indistinguishability for quantum
  particles: spin, statistics and the geometric phase}, Proc. Roy. Soc. London
  Ser. A \textbf{453} (1997), no.~1963, 1771--1790.

\bibitem{BKJ2014}
Bou Khuzam, M.N. and Johnson, M.J. \emph{On the conjectures regarding the $4$-point Atiyah determinant}. 
  SIGMA, Symmetry Integrability Geom. Methods Appl. \textbf{10} (2014), Paper 070, 9 pp.

\bibitem{Dokovic2002a}
{\Dbar}okovi{\'c}, D. \emph{Proof of Atiyah's conjecture for two special types of configurations}
Electron. J. Linear Algebra \textbf{9} (2002), 132–137.

\bibitem{Dokovic2002b}
{\Dbar}okovi{\'c}, D. \emph{Verification of Atiyah's conjecture for some nonplanar configurations with dihedral symmetry}
Publ. Inst. Math. (Beograd) (N.S.) \textbf{72}(86) (2002), 23–28.

\bibitem{Eastwood-Norbury2001}
Eastwood, M. and Norbury, P. \emph{A proof of Atiyah's conjecture on configurations of four points in Euclidean 
three-space}, Geom. Topol. \textbf{5} (2001), 885-893 (electronic).

\bibitem{Malkoun2014}
Malkoun, J. \emph{Configurations of points and the symplectic Berry-Robbins problem}, SIGMA 
Symmetry Integrability Geom. Methods Appl. \textbf{10} (2014), Paper 112, 6pp.

\bibitem{Malkoun2019}
Malkoun, J. \emph{Root systems and the Atiyah-Sutcliffe problem}, J. Math. Phys. \textbf{60}, No. 10, 101702, 6 p. (2019).

\bibitem{Malkoun2020}
Malkoun, J. \emph{Weights, Weyl-equivariant maps and a rank conjecture}, Exp. Math., published online on Jan. 27, 2020.

\bibitem{Svrtan2014}
Svrtan, D. \emph{http://www.emis.de/journals/SLC/wpapers/s73vortrag/svrtan.pdf}.

\end{thebibliography}
\end{document}